%% file: main.tex
\documentclass{article}

    \usepackage[preprint, nonatbib]{neurips_2020_tda}

\usepackage[utf8]{inputenc} 
\usepackage[T1]{fontenc}    
\usepackage{amsfonts}       
\usepackage{booktabs}       
\usepackage{hyperref}       
\usepackage{url}            
\usepackage{microtype}      
\usepackage{nicefrac}       
\usepackage{paralist}       
\usepackage{siunitx}        
\usepackage{tikz-cd}
\usepackage{amsmath}
\usepackage{amsthm}
\usepackage{bbm}
\usepackage{bm}
\usepackage{amssymb}
\usepackage{wrapfig}
\newtheorem{definition}{Definition}

\newcommand{\abs}[1]{\ensuremath{\left|{#1}\right|}}
\DeclareMathOperator{\diag}{diag}

\DeclareMathOperator{\B}{{\bf B}}
\DeclareMathOperator{\Lpn}{{\bf L}}
\DeclareMathOperator{\D}{{\bf D}}
\DeclareMathOperator{\I}{{\bf I}}
\DeclareMathOperator{\A}{{\bf A}}

\newcommand{\one}{\mathbbm{1}}

\title{
  Simplicial 2-Complex Convolutional Neural Networks
}

\author{%
  Eric Bunch \\
  Machine Learning Research \& Innovation Team \\
  American Family Insurance \\
  Madison, WI 53783 \\
  \texttt{ebunch@amfam.com} \\
  \And
  Qian You \\
  Machine Learning Research \& Innovation Team \\
  American Family Insurance \\
  Madison, WI 53783 \\
  \texttt{qyou@amfam.com} \\
  \And
  Glenn Fung \\
  Machine Learning Research \& Innovation Team \\
  American Family Insurance \\
  Madison, WI 53783 \\
  \texttt{gfung@amfam.com} \\
  \And
  Vikas Singh \\
  University of Wisconsin, Madison\\
  Madison, WI 53706\\
  \texttt{vsingh@biostat.wisc.edu}
}

\begin{document}

\maketitle

\input{abstract}

\input{introduction_and_related_work}
\input{sccnns}

\input{experiments}



\bibliographystyle{abbrv}
\bibliography{bibliography}

\end{document}

%% file: abstract.tex
\begin{abstract}
Recently, neural network architectures have been developed to accommodate when the data has the structure of a graph or, more generally, a hypergraph. While useful, graph structures can be potentially limiting. Hypergraph structures in general do not account for higher order relations between their hyperedges. Simplicial complexes offer a middle ground, with a rich theory to draw on. We develop a convolutional neural network layer on simplicial 2-complexes.

\end{abstract}

%% file: introduction_and_related_work.tex
\section{Introduction and related work}

In the last decade, neural networks have gained much popularity in the field of machine learning, in large part due to their modularity, flexibility, and ability to be successfully applied to many different types of problems. Recently there has been extensive work to formulate neural networks in settings where the data comes structured as a graph \cite{patchy_san,graphconv,diffusion_gcnn,the_gnn_model,spectral_gcnn,verma2018graph,sortpool,graphconv}, or a hypergraph \cite{hypergraph_conv,hypergraph_attention}. 

Graph structures in neural networks have proven extremely useful; however they are limited to pairwise connections of nodes. One drawback of the framework of hypergraphs is that they, in general, do not account for higher order relations among the edges of a hypergraph. Simplicial complexes offer a middle ground between graphs and hypergraphs, and have a rich theory to draw from. This document aims to begin to make the case for using the simplicial complex structure in machine learning by introducing a convolution layer on simplicial 2-complexes, akin to the popular convolution layer on graphs defined in \cite{graphconv}.
 

\subsection{Simplicial complexes}

\begin{definition}
A \emph{simplicial complex} is a collection of subsets $S$ of a finite set $X$ that is closed under taking nonempty subsets. 
\end{definition}


The finite sets comprising $S$ are called the \emph{faces} of $S$, and the dimension of a face $s$ is defined as $\text{dim}(s) = |s| - 1$. The zero-dimensional faces of $S$ will be referred to as the \emph{vertices} of $S$, and will be denoted $v_i$. Denote by $S_k$ the collection of all $k$-dimensional faces of $S$. We will write an element of $S_k$ as $s = (v_0, v_1, \dots, v_{k-1})$, where the $v_i$'s are the vertices comprising the face $s$. We say that the \emph{dimension of $S$} is the maximum dimension of all its faces. We can see that a simplicial complex of dimension one is equivalent to a graph.

Let $S$ be a simplicial complex, and assume that we have a specific enumeration of the vertices of of $S$: $v_0, \dots, v_n$. A \emph{$k$-chain} of $S$ is defined to be a finite formal sum $\sum_{i=0}^N a_is_i$, where $a_i \in \mathbb{R}$, and $s_i \in S_k$ such that $s_i = (v_{i_0}, v_{i_1}, \dots, v_{i_{k-1}})$ where $i_0 < i_1 < \cdots < i_{k-1}$. We will denote collection of $k$-chains of $S$ by $C_k$. Each $C_k$ forms a finite dimensional vector space over $\mathbb{R}$, where a choice of basis will be taken to be the collection $S_k$ of $k$-dimensional faces of $S$. If the collection $S_k$ is empty, we take $C_k$ to be the zero vector space. We also define $C_0 = 0$. For each $k = 1, 2, \dots$, there is a linear map $\partial_k: C_k \rightarrow C_{k-1}$ defined by

\begin{equation}
    \partial_k(v_0, \dots, v_{k-1}) 
    = 
    \sum_{i=0}^{k-1}(-1)^i(v_0, \dots, \widehat{v_i}, \dots, v_{k-1})
\end{equation}

\noindent where $(v_0, \dots, \widehat{v_i}, \dots, v_{k-1}) = (v_0, \dots, v_{i-1}, v_{i+1}, \dots, v_{k-1}) \in C_{k-1}$. We also define $\partial_0:C_1 \rightarrow C_0$ to be the zero map.






Given the bases described, the linear maps $\partial_k$ can be written as $ n_{k-1} \times n_k$ matrices, which we will denote by $\B_k$, following \cite{normalized_laplacian}. For each $\partial_k$, there exists a map $\partial_k^*: C_k \rightarrow C_{k+1}$, which is the adjoint of $\partial_k$. The matrix representation of $\partial_k^*$ is $\B_k^*$. 
Then we define the \emph{$k^{th}$ Laplacian} of $S$ to be the $n_k \times n_k$ matrix

\begin{equation}
    \Lpn_k = \B_k^* \B_k + \B_{k+1} \B_{k+1}^*.
\end{equation}

\noindent We can see that since $\partial_0$ is the zero map, that $\Lpn_0 = \B_1\B_1^*$. In the case when $S$ is a one dimensional complex (i.e. a graph), $\Lpn_0$ coincides with the usual definition of the graph Laplacian, which is typically defined as $L = D - A$, where $D$ and $A$ are the diagonal degree matrix and the adjacency matrix respectively of the graph.

As is the case with graph Laplacians, the $\Lpn_k$ describe how to appropriately propagate a signal on $S$. We will use these $\Lpn_k$ as well as the $\B_k$ and $\B_k^*$ construct a convolutional neural net that generalizes graph neural nets.

%% file: sccnns.tex
\section{Simplicial 2-complex convolution layer}

Let $S$ be a simplicial complex of dimension 2. For each $k = 0, 1, 2$, we write $n_k = |S_k|$, and we assume we have matrices $X_k$ of dimension $n_k \times F_k$. These will be our \emph{feature matrices}, and $F_k$ is the feature dimension. Intuitively, we think of $X_k$ as describing an $F_k$-dimensional signal on the dimension $k$ faces of $S$. One benefit of this approach is that the feature dimension, $F_k$, can vary with $k$, allowing us to encode multimodal features. 

One crucial step in this setup is properly normalizing the $\Lpn_k$. There are a variety of ways to approach this, either through normalizing the $\Lpn_k$ directly \cite{normalized_laplacian}, to use a theory of weightings on the simplicial complex itself \cite{spectra_lpn_ops}, or to appeal to a more abstract setting \cite{cellular_sheaves}. We choose here to follow \cite{normalized_laplacian} due in part to ease of implementation, and in part to the impracticality of structuring a typical data set as a high-dimensional simplicial complex.

\subsection{Normalization}

We follow \cite{normalized_laplacian}, Def. 3.3 and define the matrices

\begin{equation*}
\begin{split}
    \widetilde{\D}_2 &= \max(\diag(\abs{\B_1}\one), \I) \\
    \widetilde{\D}_3 &= \I 
    \end{split}
    \quad\quad\quad
    \begin{split}
    \D_1 &= 2\diag(\abs{\B_1}\D_2\one) \\
    \D_2 &= \max(\diag(\abs{\B_2}\one), \I) \\
    \D_3 &= \frac{1}{3}\I \\
    \end{split}
    \quad\quad\quad
    \begin{split}
    \D_4 &= \I \\
    \D_5 &= \diag(\abs{\B_2}\one) \\
\end{split}
\end{equation*}

We also follow \cite{normalized_laplacian} and define

\begin{align*}
    \begin{split}
        \widetilde{\Lpn}^u_0 
        &= 
        \B_1\widetilde{\D}_3\B_1^*\widetilde{\D}_2^{-1} \\
        \widetilde{\Lpn}^u_1 
        &= 
        \D_2\B_1^*\D_1^{-1}\B_1 \\
    \end{split}
    \begin{split}
        \widetilde{\Lpn}_1^d
        &=
        \B_2\D_3\B_2^*\D_2^{-1} \\
        \widetilde{\Lpn}^d_2 
        &= 
        \D_4\B_2^*\D_5^{-1}\B_2   
    \end{split}
\end{align*}

The $\widetilde{\Lpn}_i^u$ and $\widetilde{\Lpn}_i^d$ are called the \emph{normalized up Laplacians} and \emph{normalized down Laplacians} respectively. Define the following normalized adjacency matrices

\begin{align*}
    \begin{split}
        \A_0^u 
        &=
        \widetilde{\D}_2 - \widetilde{\Lpn}_0^u\widetilde{\D}_2
        \\
        \A_1^u 
        &= 
        \D_2 - \widetilde{\Lpn}_1^u\D_2
        \\
    \end{split}
    \begin{split}
        \A_1^d 
        &= 
        \D_2^{-1} - \D_2^{-1}\widetilde{\Lpn}_1^d
        \\
        \A_2^d 
        &= 
        \D_4^{-1} - \D_4^{-1}\widetilde{\Lpn}_2^d.
        \\
    \end{split}
\end{align*}

Finally, we define 

\begin{align*}
    \begin{split}
        \widetilde{\A}_0^u &= (\A_0^u + \I)(\widetilde{\D}_2 + \I)^{-1} \\
        \widetilde{\A}_1^u & = (\A_1^u + \I)(\D_2 + \I)^{-1} \\
    \end{split}
    \begin{split}
        \widetilde{\A}_1^d &= (\D_2 + \I)(\A_1^d + \I) \\
        \widetilde{\A}_2^d &= (\D_4 + \I)(\A_2^d + \I).
    \end{split}
\end{align*}

\noindent The $\A_i^{\alpha}$ are meant to generalize the adjacency matrix of a graph, and the $\widetilde{\A}_i^{\alpha}$ are meant to play the role of normalized adjacency matrices with added self-loops, as used in \cite{graphconv}.


\subsection{Convolution layer}
In this section, we define a simplicial complex convolutional neural network layer. Code for this implementation, as well as the experiment in \ref{sec:experiments} can be found on github \footnote{\href{https://github.com/AmFamMLTeam/simplicial-2-complex-cnns}{https://github.com/AmFamMLTeam/simplicial-2-complex-cnns}}.

We define a convolutional neural network layer using the simplicial complex structure in the following way. Assume we have hidden features at level $h$: $\{ X_k^{(h)} \}_{k=0}^2$, where $X_k^{(h)}$ is a matrix of dimension $n_k \times F_k^{(h)}$, and $X_k^{(0)} = X_k$. Then $X_k^{(h+1)}$ can be obtained by

\begin{align}
    X_0^{(h+1)} &= \sigma\left( 
        \widetilde{\A}_0^{u}X_0^{(h)}W_{0,0}^{(h)}
        \mathbin\Vert
        \D_1^{-1}\B_1X_1^{(h)}W_{1,0}^{(h)}
    \right) \label{eqn:X0h} \\ 
    X_1^{(h+1)} &= \sigma\left(
        \D_2\B_1^*\D_1^{-1}X_0^{(h)}W_{0, 1}^{(h)}
        \mathbin\Vert
        (\widetilde{\A}_1^d + \widetilde{\A}_1^u)X_1^{(h)}W_{1, 1}^{(h)}
        \mathbin\Vert
        \B_2\D_3X_{2}^{(h)}W_{2, 1}^{(h)}
    \right) \\ 
    X_2^{(h+1)} &= \sigma\left(
        \D_4\B_2^*\D_5^{-1}X_1^{(h)}W_{1, 2}^{(h)}
        \mathbin\Vert
        \widetilde{\A}_2^dX_2^{(h)}W_{2, 2}^{(h)}
    \right)
\end{align}


\noindent where $\mathbin\Vert$ denotes horizontal concatenation, $\sigma$ is an activation function, and the $W_{i, j}^{(h)}$'s are matrices of learnable weights of dimensions $F_{i}^{(h)} \times F_{j}^{(h)}$. One can see that if we have a simplicial complex of dimension one (i.e. a graph), and the features $X_1$ are all zeros, then the hidden layer output of the simplicial complex convolutional layer in Eq. \ref{eqn:X0h} becomes

\begin{equation}
    X_0^{(h+1)} = \sigma \left( \widetilde{\A}_0^u X_0^{(h)} W_{0, 0}^{(h)} \right).
\end{equation}

\noindent This coincides with the definition of a convolution layer in a graph convolutional network \cite{graphconv}.





%% file: experiments.tex
\section{Experiments}\label{sec:experiments}

Once we have all the $X_k^{(h)}$, for $k = 0, 1, 2$, and $h = 0, \dots, H$, we horizontally concatenate the features to form $X_k^{[0:H]} = X_k^{(0)} \mathbin\Vert X_k^{(1)} \mathbin\Vert \cdots \mathbin\Vert X_k^{(H)}$ to form an $n_k \times F_H$ matrix, where $F_H = \sum_{h=0}^H F_k^{(h)}$. Then we have 3 feature matrices $\{ X_k^{[0:H]} \}_{k=0}^2$. These can then be used for further downstream learning tasks, e.g. classifying faces of a certain dimension (akin to node or edge classification in graph learning), joint tasks, or classifying a simplicial complex itself when it is among a collection of simplicial complexes. The latter setting in general requires ways of dealing with simplicial complexes with differing structures. 

\subsection{MNIST classification}

We give a brief explanation of how an image can be encoded as a simplicial 2-complex with feature matrices, and employ simplicial complex convolution layers in a network to classify handwritten digits in the MNIST data set. We detail a proof of concept for how these simplicial complex convolution layers can be used to augment traditional convolution layers through an ablation study.

\begin{wrapfigure}{r}{0.4\textwidth}
\centering
\includegraphics[width=0.30\textwidth]{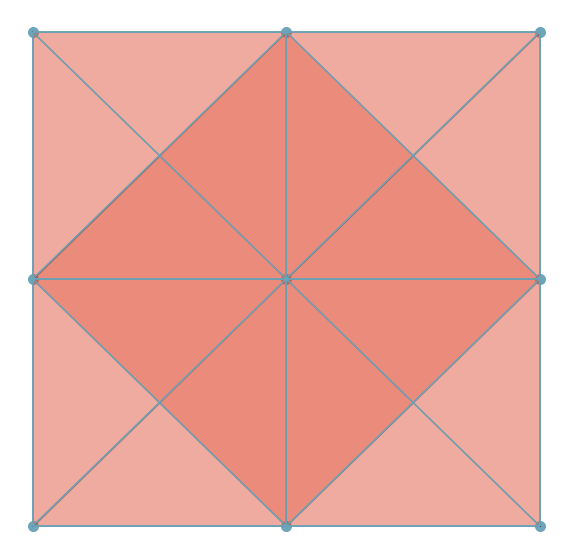}
\caption{Depiction of the type of 2-simplicial complex structure constructed for the experiments described.}
\label{fig:sc_grid}
\end{wrapfigure}

Given an image of size $h \times w$ pixels, and a choice of \textit{kernel size} $k$ and \textit{step size} $s$ (the same idea from traditional CNNs), we form a grid of zero faces that is $\left(\lfloor \frac{h - k}{s} \rfloor + 1 \right) \times \left( \lfloor \frac{w - k}{s} \rfloor + 1 \right)$. We then add 1-faces to connect any zero faces that are horizontally, vertically, or diagonally adjacent. We then add two faces wherever possible; that is, a triplet of zero faces holds a two face whenever they are all pairwise connected by a one face. Figure \ref{fig:sc_grid} gives an example depiction of such a simplicial complex. The feature vector attached to a zero face is the $k^2 \times 1$ vector resulting from flattening the $k \times k$ square of pixel values corresponding to that zero face. These are collected into an $n_0 \times k^2$ matrix $X_0$, with rows indexed by zero faces, and columns indexed by feature values. To both one faces and two faces, the features attached are $1 \times 1$ vectors with the value 1.

The ablation study takes as the base a neural net that has one traditional convolution layer, followed by a fully connected network with one hidden layer. Next, a single graph convolution layer is prepended to this network to form the first comparison. The graph is taken to be the one formed by taking only the zero and one faces of the simplicial complex described above. For the final comparison, a single simplicial complex convolution layer is prepended to the initial architecture described.

We describe first a particular 1D convolution layer. For a feature matrix $X$ of dimension $n \times F$, we flatten $X$ to a matrix of dimension $1 \times nF$. We denote by ${\bf CONV1D}_F$ a 1D convolution layer with kernel size of $F$ and step size $F$. We write ${\bf CONV1D}$ when the input feature size is understood. We can see that in the setting described above, ${\bf CONV1D}_{k^2}$ applied to a flattened $X_0$ (of size $1 \times n_0k^2$) is equivalent to a 2D convolution layer with kernel size $k$ and step size $s$.

We denote by ${\bf FC}$ a fully connected network with one hidden layer of size 32, and output layer of size 10 (for the 10 classes of MNIST). We denote by ${\bf GCONV}$ a graph convolutional layer with embedding dimension 32. We denote by ${\bf SCCONV}$ a simplicial complex convolutional layer with embedding dimensions of 32 for each face dimension. In the ablation study, the combination ${\bf SCCONV-CONV1D}$ applies a 1D convolution layer to each of the three feature matrices; one for each face dimension.

Each model was trained (tested) on a stratified sample of 1000 images from the designated training (testing) portion of the MNIST data set. Each model was trained for 300 epochs with a batch size of 8, learning rate of $10^{-4}$, dropout rate of 10\% on the fully connected layers, and batch normalization. The kernel size and step size are both taken to be 4 for all models. The metrics reported in Table \ref{tab:ablation_study} are the average accuracy, with the standard deviation denoted in the usual fashion, for five iterations of the experiment setup described above. On each experiment, the test and train set are sampled anew.

\begin{table}[ht]
    \centering
    \begin{tabular}{c c}
        \toprule
             {\bf Model} & {\bf acc. $\pm$ std.} \\
             \midrule
             ${\bf CONV1D-FC}$ & 87.40 $\pm$ 0.70 \\
             \addlinespace[0.5ex]
             ${\bf GCONV-CONV1D-FC}$ & 87.42 $\pm$ 0.59 \\
             \addlinespace[0.5ex]
             ${\bf SCCONV-CONV1D-FC}$ & {\bf 91.10} $\pm$ 0.40 \\
        \bottomrule
        \addlinespace[2ex]
    \end{tabular}
    \caption{Ablation study results. Statistics are calculated over five experiments.}
    \label{tab:ablation_study}
\end{table}

We note that these metrics are not state of the art, and there is still a great deal to still investigate in this direction. The results of this ablation study are suggestive enough to warrant further study of the simplicial complex convolution layer in general. 
